\ifx\selectfont\undefined
\documentstyle[12pt]{article}
\else
\documentstyle[12pt,oldlfont]{article}
\fi

\newtheorem{thm}{Theorem}
\newtheorem{prop}{Proposition}

\newcommand{\rem}{\noindent{\bf Remark{\ \ }}}
\newcommand{\proof}{{\noindent\bf Proof{\ \ }}}
\newcommand{\qed}{\mbox{}\hfill\(\Box\)\bigskip}

\newcommand{\Rn}[1]{\mbox{{\it I\kern -0.25emR}$^{\,{#1}}$}}
\newcommand{\NN}{\mbox{{\it I\kern -0.25emN}}}
\newcommand{\al}{\alpha}

\newcommand{\De}{\Delta}
\newcommand{\ep}{\varepsilon}

\newcommand{\spn}{\mathop{\overline{\rm span}}\,}
\newcommand{\lust}{\mathop{\rm \, l.u.st.}}

\newcommand{\restr}[2]{#1 \,\vert\sb {#2}}
\title{Erratum to: ``Banach Spaces without\\
Local Unconditional Structure''}
\author{Ryszard A.~Komorowski
        \and  Nicole Tomczak-Jaegermann
}
\newcommand\address{\noindent\leavevmode%
{\small
Institute of Mathematics,
Technical University,\\
Wroc\l aw, Poland,\\
{\small\tt%
komo@graf.im.pwr.wroc.pl}\\[.5cm]
\noindent
Department of Mathematics,
University of Alberta,\\
Edmonton, Alberta, Canada T6G 2G1,\\
{\small\tt%
  ntomczak@approx.math.ualberta.ca} }}
\date{7 July, 1996}
\begin{document}
\maketitle

\begin{abstract}
  This note contains a corrected proof of the main result (which
  remains unchanged) from [K-T].  It was recently observed that an
  argument in a basic technical criterium has a gap.
\end{abstract}

The first named author recently observed that the paper by Komorowski
and Tomczak [K-T], contains a gap.  Here we shall correct the
formulation of a general result (Theorem 2.1 in [K-T]) and then we shall
complete the proof of the main application to the homogeneous Banach
space problem (Theorem 3.1 in [K-T]), that is now neccessary in view of
modifications in the aforementioned general statement.
Corollaries 4.3, 4.4 and 4.5  do not require any changes.

This note is not self-contained. However, in order to make it possible
for the readers a little familiar with [K-T], to read it without
constantly referring to the original paper, we decided not to cut the
arguments to bare minimum, and occasionally give some additional
background information.

\bigskip

The gap is in the proof of the fundamental criterium from [K-T]
(Proposition 1.1).  One deals there with a Banach space $Y$ which has a
2-dimensional Schauder decomposition $\{Z_{k}\}_k$. Assuming that $Y$
has cotype $q < \infty$ and $\lust$, one constructs a block diagonal
operator $T$, which was claimed to be bounded on subspaces $Y_K$ for
which the decomposition $\{Z_{k}\}_{k \in K}$ was unconditional.  In
fact, we can estimate this norm only if the subspace is complemented by
a natural projection. To be more precise, for a subset $K \subset \NN$
we define $Q_K: Y \to Y_K$ by $Q_K (\sum_k z_k)= \sum_{k \in K} z_k$.
Then the correct formulation of Proposition 1.1 states, in condition
(ii), that if the decomposition $\{Z_{k}\}_{ k \in K}$ of $Y_K$ is
$C$-unconditional and $\|Q_K\|\le C'$, then
\begin{equation}
\|\restr{T}{Y_K}\colon Y_K \to Y_K \| \le  C^2 C'
       \psi  \lust(Y),
  \label{borz_ii}
\end{equation}
where $\psi$ is the same function as in [K-T].

The proof uses the same operator $T$ as in [K-T], and the estimate is
completely clear from the original proof, once we notice that letting
$P_k\colon Y \to Z_k$ to denote the natural projection onto $Z_k$, for
$k=1, 2, \ldots$, we trivially have $P_k = P_k Q_{K}$ for every $k \in
K$, and thus, for arbitrary signs $\ep_k = \pm 1$ we have
$$
\Bigl\|\sum_{k\in {K}} \ep_k  P_k : Y \to Y_{K} \Bigr\| =
\Bigl\|\sum_{k\in {K}} \ep_k  P_k Q_{K}\Bigr\|\le C\,C'.
$$ 

\rem 
Note that if the decomposition $\{Z_{k}\}_k$ of $Y$ is unconditional,
the natural projections $Q_K$ are automatically bounded for all $K
\subset \NN$. In this case the proof of Proposition 1.1 in [K-T] was
correct, so Corollaries 4.3, 4.4 and 4.5  do not require any changes.

\bigskip

The remaining part of this note heavily depends on  technical  notation 
from the end of Section 1 of [K-T]. For the reader convenience
we copy it here.

If $\Delta = \{A_{ m}\}_m$ and $\Delta' = \{A_{ m}'\}_m$ are
partitions of $\NN$, we say that $\Delta \succ \Delta'$, if there
exists a partition ${\cal J} (\Delta', \Delta)= \{J_{ m}\}_m$ of $\NN$
such that
\begin{equation}
\label{suc_part_a}
\min J_m < \min J_{m+1}
\qquad \mbox{and} \qquad
A_{m}' = \bigcup_{j \in J_m} A_{j}
\qquad \mbox{for\ } m=1, 2,\ldots.
\end{equation}
In such a situation, for $m=1, 2, \ldots$,
${\cal K}(A_m', \Delta)$  denotes the family
\begin{equation}
\label{suc_part_b}
{\cal K}(A_m', \Delta) = \{ K \subset A_{m}'\, \mid\,
|K \cap  A_{j}|=1 \ \mbox{for\ } j \in J_m\}.
\end{equation}

Finally, if  $\Delta_i = \{A_{i, m}\}_m$,
for $i= 1, 2,\ldots$,  is a sequence of partitions of $\NN$,
with $\Delta_{1} \succ \ldots \succ \Delta_{i} \succ \ldots$,
we set, for $m=1, 2, \ldots$ and $i= 2,3, \ldots$
\begin{equation}
{\cal K}_{i, m}= {\cal K} (A_{i, m}, \Delta_{i-1}).
  \label{k_im}
\end{equation}

Now Theorem 2.1 from [K-T] requires an obvious extra assumption of
boundness of natural projections. Although the statement of this
theorem is rather lengthy, we recall it here (in a corrected form) for
sake of future references.

\begin{thm}\mbox{\bf $\!\!\!\star$2.1:}
Let $X= F_{1}\oplus\ldots\oplus F_{4}$
be a direct sum of  Banach spaces
of cotype $r$, for some $r < \infty$,
and let $\{f_{i,l}\}_{l}$ be a  normalized
monotone Schauder basis
in $F_i$, for $i=1, \ldots, 4$.
Let $\Delta_{1} \succ \ldots \succ \Delta_{4}$
be partitions of $\NN$,
$\Delta_i = \{A_{i, m}\}_m$  for $i=1, \ldots, 4$.
Assume that  there exist  $C, C' \ge 1$ such that
for every $K \in {\cal K}_{i, m}$  with
$i=2,3, 4$ and  $m=1, 2, \ldots$,
 the  basis
$\{f_{s,l}\}_{l \in K}$ in $\restr{F_s}{K}$ is
$C$-unconditional, and the natural projection
$R_K^{(s)}\colon F_s \to \restr{F_s}{K}$ has
the norm  $\| R_K^{(s)}\| \le C'$, 
for $s=1,  \ldots,4$;  moreover,
there is $\widetilde{C} \ge 1$ such that
for  $i=1, 2, 3$  and  $m=1, 2, \ldots$
we have
\begin{equation}
\|I\colon \restr{F_i}{A_{i,m} }\to
     \restr{F_{i+1}}{A_{i, m}} \|\le \widetilde{C}.
  \label{4_2_small}
\end{equation}
Assume finally  that one of the following
conditions is satisfied:
\begin{description}
\item[(i)]
there is a sequence $0 < \delta_m <1$ with
$\delta_m \downarrow 0$  such that
for every $i= 1, 2, 3$ and  $m=1, 2, \ldots$
and  every  $K \in {\cal K}_{i+1, m}$ we have
\begin{equation}
\|I\colon \restr{D(F_{1}\oplus\ldots\oplus F_{i})}{K}
\to \restr{F_{i+1}}{K}\|  \ge \delta_m ^{-1};
\label{4_2_large}
\end{equation}
\item[(ii)]
there is a sequence $0 < \delta_m <1$ with
$\delta_m \downarrow 0$ and
$\sum_m \delta_m ^{1/2} = \gamma <\infty$
such that
for every $i= 1, 2, 3$ and  $m=1, 2, \ldots$
and  every  $K \in {\cal K}_{i+1, m}$ we have
\begin{equation}
\|I\colon \restr{F_{i+1}}{K}
    \to \restr{F_{i}}{K}\|  \ge \delta_m ^{-1}.
\label{4_3_large}
\end{equation}
\end{description}
Then there exists a subspace $Y$ of  $X$
without local  unconditional structure,
but which still admits a Schauder basis.
\end{thm}

\proof
The only addition to the original proof is to check that the assumption
for $\|Q_K\|$ required  in  the fundamental criterium is now
satisfied. Since  spaces $Z_k$ are spanned by vectors of the form
\begin{eqnarray*}
  x_k &=& \al_{1, k} f_{1, k}+ \ldots + \al_{4, k} f_{4, k},\\
  y_k &=& \al_{1, k}' f_{1, k}+ \ldots + \al_{4, k}' f_{4, k},
\end{eqnarray*}
then  we have 
$Q_K ( \sum_k z_k) = \sum_{s=1}^4 R_K^{(s)} \Bigl(\sum_{k} 
   (t_k \al_{s, k} + t'_k \al_{s, k}') f_{s, k}\Bigr)$,
for  $z =  \sum_k z_k = \sum_k (t_k x_k + t'_k y_k)$,
thus $\|Q_K\|\le 4 \max_s \|R_K^{(s)}\|$.
\qed

\bigskip

The remaining part of this note contains the proof of Theorem 3.1 from
[K-T], which says:

\begin{thm}\mbox{\bf $\!\!\!$3.1:}
Let $X$ be a Banach space with an unconditional basis
and of cotype $r$, for some $r < \infty$.
If $X$ does not contain a subspace isomorphic  to $l_2$ then
there exists a subspace $Y$ of  $X$ without local unconditional
structure, which admits a Schauder basis.
\end{thm}

Let $\{e_{l}\}_{l}$ be a 1-unconditional basis in $X$.  Passing to a
subsequence we may assume without loss of generality that
$\{e_{l}\}_{l}$ generates a spreading model $\{u_{l}\}_{l}$.  We
distinguish two main cases.
\begin{description}
\item[(I)] either no sequence of  disjointly supported
  blocks of $\{e_{l}\}_{l}$ of length $\le 3$ with constant
  coefficients is equivalent to $\{u_{l}\}_{l}$,
\item[(II)] or the basis  $\{e_{l}\}_{l}$ is 1-subsymmetric.
\end{description}

If Case I does not hold, then there exist disjointly supported blocks
$\{w_{l}\}_{l}$ equivalent to $\{u_{l}\}_{l}$; so by renorming we get
into Case II for $X_1 = \spn \{w_{l}\}_{l}$.  Case I is new, compared
with [K-T].  Case II has the same proof as in [K-T].  The additional
information that the basis is subsymmetric implies that equally
distributed vectors spanning $\ell_2^n$'s are given by constant
coefficient blocks (see Lemma 3.2 in [K-T] and the remark after). Then
the required complementation follows from the classical fact that in a
subsymmetric space subspaces spanned by constant coefficient
consecutive blocks are well complemented.

\medskip

Case I is the content of the following proposition proved by
B.~Maurey and N.~Tomczak in May 1993 (unpublished).

\begin{prop}\mbox{\bf $\!\!\!$I:}
  Let $X$ be a Banach space of cotype $r$, for some $r < \infty$.
  Assume that $X$ has a 1-unconditional basis $\{e_{l}\}_{l}$ which
  generates a spreading model $\{u_{l}\}_{l}$. If no disjointly
  supported blocks of $\{e_{l}\}_{l}$ of length $\le 3$ with constant
  coefficients  are  equivalent to $\{u_{l}\}_{l}$, then there exists a
  subspace $Y$ of $X$ without local unconditional structure, which
  admits a 2-dimensional unconditional decomposition
\end{prop}

\proof
Denote the span  of  $\{u_l\}_l$ by $E$.
First observe, by a compactness argument, that if no subsequence
of $\{e_l\}_l$ dominates (resp. is dominated  by) $\{u_l\}_l$
and if $\{A_m\}_m$ is a partition of $\NN$ into finite sets
then for every $D$ there exists $M$ such that whenever
$K $ is a subset of $\NN$ such that $|K \cap A_m|=1$
for $ 1 \le m \le M$ 
then  $\|I: \restr{X}{K} \to \restr{E}{K}\| \ge D$
(resp. 
$\|I: \restr{E}{K} \to \restr{X}{K}\| \ge D$).  

\smallskip
By passing to subsequences of $\{e_l\}_l$ we may assume that 
\begin{description}
\item[(Ia)] either
$\|I \colon \restr{E}{L} \to \restr{X}{L}\| = \infty$ for all 
infinite subsequences $L$ of $\NN$,
\item[(Ib)] or $\|I \colon E \to X \| < \infty$.
\end{description}

Partition the whole basis $\{e_l\}_l$ into four infinite sets
$\{e_{i,l}\}_l$  for $i=1, \ldots, 4$. Of course $E$ is the spreading
model  for each  $\{e_{i,l}\}_l$ ($i=1, \ldots, 4$).
In  the proof we shall use a shorthand notation
$\sim$ to denote that two sequences are 2-equivalent.

\medskip\noindent 
{\bf Case (Ia):\ } 
Set $f_{1, l} = e_{1,l}$ for all $l$ and let $\Delta_1 = \{A_{1,
  m}\}_m$ consists of singletons. Let $F_1 = \spn [f_{1, l}]_l$.  Let
$\Delta_2 = \{A_{2, m}\}_m$ be a partition of $\NN$ into successive
intervals such that $\|I: \restr{E}{A_{2, m}} \to
\restr{F_1}{A_{2,m}}\| \ge 2^{2m} $ for all $m$. Define $f_{2, l}$
for $l \in A_{2, m}$ by induction in $m$.  For $l \in A_{2, 1}$ set
$f_{2, l}= e_{2, s_1 +l}$, where $s_1$ is so large that $\{f_{2,
  l}\}_{l \in A_{2, 1}} \sim \{u_l\}_{l=1}^{|A_{2, 1}|}$.  Then having
defined $f_{2, l}$ for all $l \in A_{2, m}$ set, for $l \in A_{2,
  m+1}$, $f_{2, l}= e_{2, s_{m+1} +l}$, where $s_{m+1} > s_m + |A_{2,
  m}|$ is so large that $\{f_{2, l}\}_{l \in A_{2, m+1}} \sim
\{u_l\}_{l=1}^{|A_{2, m+1}|}$.  In particular, $\{f_{2, l}\}_l$ is a
subsequence of $\{e_{2,l}\}_l$.  Let $F_2 = \spn [f_{2, l}]_l$.

By our initial remark we can define by induction a sequence
$0 = M_0 < M_1 < \ldots < M_j < \ldots $ such that for every 
$j=1, 2, \ldots  $ whenever $K \subset \NN$ satisfies
$|K \cap A_{2, m}|=1$ for $ M_{j-1}< m \le M_j$, 
then $\|I: \restr{E}{K} \to \restr{F_2}{K}\| \ge 2^{2j}$.
Then set $A_{3, j} = \bigcup_{m= M_{j-1}+1}^ {M_j} A_{2, m}$
for $j= 1, 2, \ldots$.
Define $\{f_{3, l}\}_l$ as a subsequence of $\{e_{3, l}\}_l$
such that 
$\{f_{3, l}\}_{l \in A_{3, j}} \sim \{u_l\}_{l=1}^{|A_{3, j}|}$,
and  set $F_3 = \spn [f_{3, l}]_l$. 

Finally, using the same construction once more for the partition
$\Delta_3= \{A_{3, j}\}_j$, find a partition $\Delta_4 = \{A_{4, n}\}$
with $\Delta_4 \prec \Delta_3$ such that for every $K \in {\cal K}_{4,
  n}$ we have $\|I: \restr{E}{K} \to \restr{F_3}{K}\| \ge 2^{2n}$ (for
$n=1, 2, \ldots$).  Then define $\{f_{4, l}\}_l$ as a subsequence of
$\{e_{4, l}\}_l$ such that $\{f_{4, l}\}_{l \in A_{4, n}} \sim
\{u_l\}_{l=1}^{|A_{4, n}|}$ for $n=1, 2, \ldots$; and by $F_4 = \spn
[f_{4, l}]_l$ denote the corresponding subspace.

To check  conditions  (\ref{4_3_large}), it is enough 
to observe that  for $i=1, 2, 3$ and $m=1, 2, \ldots$, and   every  
$K \in {\cal K}_{i+1, m}$ we have $K \subset A_{i+1, m}$;
and hence $\{f_{i+1, l}\}_{l \in K} \sim \{u_l\}_{l=1}^{|K|}$,
in other words, $\restr{F_{i+1}}{K} \sim \restr{E}{K}$.
Therefore,  by our construction,
$$
\|I\colon \restr{F_{i+1}}{K} \to \restr{F_{i}}{K}\| 
\ge (1/2) \|I\colon \restr{E}{K} \to \restr{F_{i}}{K}\| 
\ge 2^{2m-1}= \delta_m ^{-1}.
$$

To check  (\ref{4_2_small}), note that 
for  $i=1, 2, 3$  and  $m=1, 2, \ldots$ and 
every $A_{i,m}$  one has
$\restr{F_{i}}{A_{i, m}} \sim \restr{F_{i+1}}{A_{i, m}} 
\sim \restr{E}{A_{i, m}}$.

Finally, all the $F_s$'s are spanned by subsequences of the 
1-unconditional basis, hence  the projections $R_K^{(s)}$
have norms equal to 1.
This shows case (Ia).

\medskip\noindent 
{\bf Case (Ib):\ } 
Set $f_{1, l} = e_{1,l}$ for all $l$ and let $\{A_{1, m}\}_m$ be
singletons. Let $F_1 = \spn [f_{1, l}]_l$.  Since $\|I\colon E \to F_1
\| \le \|I: E \to X \| < \infty$, we must have $\|I: F_1 \to E \| =
\infty$.  Let $\Delta_2 = \{A_{2, m}\}_m$ be a partition of $\NN$ into
successive intervals such that $\|I \colon \restr{F_1}{A_{2, m}} \to
\restr{E}{A_{2,m}}\| \ge 2^{2m} $ for all $m$.  Define the sequence
$\{f_{2, l}\}_l$ as a subsequence of $\{e_{2, l}\}_l$ in the same way as in
Case (Ia) and let $F_2$ be its span. Then for all $m$ we have
$\restr{F_2}{A_{2, m}} \sim \restr{E}{A_{2, m}} $, hence $\|I \colon
\restr{F_1}{A_{2, m}} \to \restr{F_2}{A_{2,m}}\| \ge 2^{2m} $.

Since 
$\|I\colon E \to F_1 \|  < \infty$, and 
$\|I\colon E \to F_2 \|  < \infty$, then also
$\|I\colon E \to D(F_1\oplus F_2) \|  < \infty$, 
where $  D(F_1\oplus F_2) = \spn [f_{1, l} +  f_{2, l}]_l$.
Thus  $\|I\colon D(F_1 \oplus F_2) \to E \|  = \infty$.
As in Case (Ia) we can find an increasing  sequence $\{M_j\}$ of
integers and a partition  $\Delta_3 = \{A_{3, j}\}_j$ of $\NN$ with
$A_{3, j} = \bigcup_{m= M_{j-1}+1}^ {M_j} A_{2, m}$ such that 
for every $K \in {\cal K}_{3, j}$ we have
$\|I\colon \restr{D(F_{1}\oplus F_{2})}{K}
\to \restr{E}{K}\|  \ge 2^{2j}$.

We then  define a subsequence $\{f_{3, l}\}_l$ of $\{e_{3, l}\}_l$,
so that its span $F_3$ satisfies 
$\restr{F_3}{A_{3, j}} \sim \restr{E}{A_{3, j}} $,  for all $j$.

Finally we define a partition $\Delta_4 = \{A_{4, n}\}_n$ with
$\Delta_4 \prec \Delta_3$ such that for every $K \in {\cal K}_{4, n}$
we have $\|I \colon \restr{ D(F_{1}\oplus F_{2}\oplus F_3)}{K} 
\to \restr{E}{K}\| \ge 2^{2n}$ (for $n=1, 2, \ldots$), and then
we define a subsequence $\{f_{4, l}\}_l$  of $\{e_{3, l}\}_l$,
so that its span $F_4$ satisfies 
$\restr{F_4}{A_{4, n}} \sim \restr{E}{A_{4, n}}$,  for all $n$.

Just as in Case (Ia), it is easy to check that (\ref{4_2_small})
and (\ref{4_2_large}) are satisfied, and that all the projections
$R_K^{(s)}$  have norms equal to 1.
\qed

\medskip
\noindent{\bf Proof of Theorem 3.1 in Case (II):}
The key argument in [K-T] was contained in Proposition 3.3. The same
statement and the proof that followed, work here, provided that we can
use Theorem $\star$2.1. This requires that the notion of a $C$-regular
pair $\{\Delta, F\}$ used in [K-T] should be replaced by a suitable
modification, which we call here a $C$-projection-regular pair.

We say that a partition $\Delta= \{A_{m}\}_m$ of $\NN$ into
consecutive intervals and a space $F$ with a normalized Schau\-der
basis $\{f_{l}\}_{l}$ form a $C$-projection-regular pair, if the
following conditions are satisfied:
\begin{description}
\item[(i)] $ \mbox{equiv} \Bigl(\restr{F }{A_{m} },\, l_2 ^{|A_{m} |}\Bigr)
  \le   C $  and the natural projection $Q_{A_m}\colon F \to \restr{F
    }{A_{m} }$ has the norm $\|Q_{A_m}\| \le C$, for $m=1, 2, \ldots$;
\item[(ii)] for every $L \in {\cal L}(\De) $, the basis $\{f_{l}\}_{l
    \in L}$ in $\restr{F }{L}$ is 1-unconditional and the natural
  projection $Q_L \colon F \to \restr{F }{L} $ has the norm $\|Q_L \|
  \le C$ (here $ {\cal L}(\Delta )$ is the family of all $ L \subset
  \NN$ such that $|L \cap A_m|=1 $ for all $m$).
\item[(iii)] for arbitrary $L, L' \in {\cal L}(\Delta ) $ one has
$\mbox{equiv} \Bigl(\restr{F }{L}, \restr{F }{L'}\Bigr) =1$.
\end{description}

The modified Proposition 3.3 now additionally assumes that all spaces
$E_i$ have subsymmetric bases, and it asserts that the constructed
pairs $\{\Delta_i, F_i\}$ are $C$-projection-regular.

\smallskip 

The use of Rademacher functions in a setting of discrete Banach
lattice gives more information for spaces with subsymmetric bases. For
example, Tzafriri's argument mentioned in [K-T] uses this approach;
this can be also tackled by a modification of inequalities from Lemma
3.2 in [K-T]. One gets that for $m \in \NN$ and $N = 2^m$, if $E$ is an
$N$-dimensional space of cotype $r < \infty$ with  a 1-subsymmetric basis
$\{e_j\}_j$, then $m$ Rademacher vectors $\{f_l\}$ of the form
\begin{equation}
  \label{radem}
 f_l = \alpha \sum_{j=1}^N \ep_l(j) e_j
\end{equation}
are $C$-equivalent to the unit vector basis in $\ell_2^m$, where $C$
depends on $r$ and the cotype $r$ constant. Moreover, $\{\ep_l(j)\}$
are values of Rademacher functions on the diadic partition of the
interval $[0, 1]$, and so they are mutually orthogonal, $\sum_j
\ep_l(j) \ep_{l'}(j) =0$ if $l \ne l'$.

\smallskip

For a given space $E$ of cotype $r$ with a subsymmetric ba\-sis
$\{e_j\}_j$, and a partition $\Delta = \{A_{m}\}_m$ of $\NN$ into
consecutive intervals, we repeate word by word the first part of the
proof of Proposition 3.3 in [K-T] to get a subspace $F = \spn [f_l]_l$
of the following form: we fix  arbitrary successive intervals of integers,
$\{J_m\}$, with $|J_m| = 2^{|A_m|}$ and  by formula (\ref{radem})
we define vectors $f_l$ by 
$$
f_l = \alpha_m \sum_{j \in J_m} \ep_l(j) e_j
\qquad \mbox{for } l \in A_m,\ m=1, 2, \ldots.
$$ 

All the conditions from [K-T] are satisfied as before, we need only to
check the boundness of the natural projections, $Q_K\colon F \to
\restr{F}{K}$ given by $Q_K (\sum_l a_l f_l) = \sum_{l\in K} a_l f_l$,
for $K \subset \NN$.

For condition (i), $Q_{A_m}$ is equal to the restriction to $F$ of the 
natural projection $P_{J_m}$ in $E$ onto  coordinates from $J_m$,
hence is obviously  bounded.

For condition (ii), let $L = \{l_m\}_m \in {\cal L}(\De) $. Then
$\{f_{l_m}\}_m$ form a subspace of $E$ spanned by constant coefficient
successive blocks.  Consider the averaging projection $P_L$ in $E$
given by the formula
$$
P_L \bigl(\sum_j a_j e_j\bigr) = \sum_m
(\al_m |J_m|)^{-1} \Bigl( \sum_{j \in J_m}\ep_{l_m}(j) a_j\Bigr) f_{l_m}.
$$ 

The orthogonality relations of sign vectors $\{\ep_l(j)\}$ imply that
for any $l \not\in L$ we have $P_L (f_l) = 0$. So the projection $Q_L$
on $F$ is equal to the restriction  of $P_L$ to $F$. Finally let us
recall (Lindenstrauss--Tzafriri, Vol I, 3.a.4) that averaging
projections in a space with a subsymmetric basis have the norm
$\|P_L\|\le 2$, and this shows  $\|Q_L\|\le 2$.

This shows a modified version of Proposition 3.3 which allows to use
Theorem $\star$2.1 in the proof of Theorem 3.1.
\qed

\rem 
Actually, even in Case II, a subspace without l.u.st. can be constructed
which has 2-dimensional unconditional decomposition.  This can be done
by a careful use of Krivine theorem on finite block representability of
$\ell_p$. Since the basis $\{e_l\}_l$ is subsymmetric, one can find $1
\le p \le \infty$ and an unconditional basic sequence $\{f_l\}_l$ in $X=
\spn [e_l]_l$ such that denoting by $\Delta = \{J_m\}_m$ the partition
of $\NN$ into successive intervals with $|J_m| = m$, we get $\{f_l\}_{l
  \in J_m} \sim \ell_p^m$ and for any $L, L' \in {\cal L}(\Delta)$, the
sequences $\{f_l\}_{l \in L}$ and $\{f_l\}_{l \in L'}$ are 1-equivalent,
but they are not equivalent to $\ell_p$.  In many cases, $\{f_l\}_l$ is
a block basis consisting of Krivine's vectors; and then, for $l \in
J_m$, the $f_l$'s are shifts of the same vector along the basis
$\{e_l\}_l$. In remaining cases, for instance if $X= \ell_q$ for some $q
\ne 2$, we let $p =2$ and use Dvoretzky theorem. After constructing a
suitable basis $\{f_l\}_l$, the construction follows lines similar as in
Case I.

We feel, however, that the construction of Case II based on [K-T] may
be of interest in other situations then Theorem 3.1 as well.
Similarly, Theorem $\star$2.1 is, as far as we know, one of very few
results of this type which does not assume the unconditionality of the
whole Schauder decomposition, but just of its subsets.

\bigskip
\noindent{\bf References:}

\noindent{\bf [K-T]\ }R.~Komorowski \& N.~Tomczak-Jaegermann:
  Banach spaces without local unconditional structure; 
{\sl Israel J.~of Math.}, 89 (1995), pp.  205--226. 

\bigskip

\address

\end{document}